\documentclass [smallextended,draft,nospthms] {svjour3} [06-May-2014]
\usepackage{undertilde}
\usepackage{amsmath,amsfonts}
\usepackage{amsthm}
\usepackage{graphicx}
\usepackage{textcomp}
\usepackage{breqn}
\usepackage{enumerate}
\theoremstyle{plain}
  \newtheorem{theorem}{Theorem}

   \newtheorem{lemma}[theorem]{Lemma}

\title{
The Logarithmic Super Divergence and Statistical Inference 
: Asymptotic Properties
}
\date{}

\author{
Avijit Maji, Abhik Ghosh and Ayanendranath Basu \\
}
\institute {Indian Statistical Institute
\at 203, B.T. Road, Kolkata-700108, India.\\
\email{avijit.maji$@$hotmail.com, abhianik@gmail.com, ayanendranath.basu@gmail.com.}
}
\begin{document}
\maketitle

\begin{abstract}
Statistical inference based on divergence measures have 
a long history.
Recently, Maji, Ghosh and Basu (2014) 
have introduced a general family of divergences 
called the
logarithmic super divergence (LSD) family. This family acts 
as a superfamily for both of the logarithmic power 
divergence (LPD) family (eg. Renyi, 1961) and the logarithmic density
power divergence (LDPD) family introduced by 
Jones et al.~(2001). 
In this paper we describe the asymptotic properties 
of the inference procedures resulting from this 
divergence in discrete models. The properties are well supported 
by real data examples.
\end{abstract}
\vspace{2 mm}

\section{Introduction}
\label{Sec_Introduction}

The 
density-based minimum divergence approach 
has long been an important parametric inference tool.
In this approach the 
closeness between the data and the model is measured by
a density-based divergence between the data density and the
model density, such as 
a $\chi^2$ type divergence or a $\phi$-divergence
(Csisz\'{a}r 1963, 1967a,b; Ali and Silvey 1966)
 or a Bregman divergence (Bregman, 1967). 
Apart from their natural appeal, most
of these 
methods are very useful for their inherent robustness properties.
A prominent member of the class of density-based divergences 
is the Pearson's $\chi^2$ (Pearson, 1900)
which started its journey 
from the very early days
of formal research
in statistics. 
From the robustness perspective, however, Beran's 1977 
work is the first useful reference in the literature
of density-based minimum divergence inference. 
In the present paper
we focus on some variants of the power divergence (PD) measure of 
Cressie and Read (1984) and the density power divergence (DPD) of Basu et al.~(1998)
and discuss various properties related to statistical inference based on a
generalized superfamilies of these measures. 
The primary divergence class 
of logarithmic super divergences (LSDs)
which is of interest to us in this paper 
has been proposed recently by Maji, Ghosh and Basu (2014). 
In the present paper we establish the theoretical asymptotic
properties of the resulting statistical procedures. 
\paragraph{}
The rest of the paper is organized as follows. Section \ref{SEC:LPD_LDPD}
describes the logarithmic power divergence family and the 
logarithmic density power divergence family. Section \ref{SEC:lsd_defn}
gives the form and the estimating equation of the 
logarithmic super divergence (LSD) family whereas Section \ref{SEC:lsd_properties}
establishes the asymptotic distribution for LSD estimator. Section \ref{SEC:est_example}
gives some illustration of the proposed estimation procedure with real data. 
Section \ref{SEC:testing} gives the asymptotic
distribution of the test statistic using LSD measure for both one sample 
and two sample problem and Section \ref{SEC:testing:example} 
provides a hypothesis testing example.
Concluding remarks are in Section \ref{Sec_Conclusion}.
\section{The Logarithmic Power Divergence (LPD) and the Logarithmic 
Density Power Divergence (LDPD) Families}\label{SEC:LPD_LDPD}
Jones et al.~(2001) described 
a class of divergence measures
which 
do not require 
any nonparametric smoothing techniques for their construction. 
This 
family of divergences is given by 
\begin{equation}\label{EQ:LDPD}
{\rm LDPD}_\beta(g,f) = \displaystyle \log \int  f^{1+\beta} - 
\left(1 + \frac{1}{\beta}\right)  
\log \int f^\beta g + \frac{1}{\beta} \log \int g^{1+\beta},
\end{equation}
$\beta\geq0$,
where LDPD stands for logarithmic density power divergence. 
For the case $\beta=0$ we define 
$${\rm LDPD}_0(g,f)=
\lim_{\beta\rightarrow0} {\rm LDPD}_\beta(g,f)=
\lim_{\lambda\rightarrow0} {\rm PD}_\lambda(g,f)={\rm PD}_0(g,f),$$
where the power divergence (PD) measure has the form
\begin{equation}\label{EQ:cressie-read}
{\rm PD}_\lambda(g,f)=\frac{1}
{\lambda(\lambda + 1)} \int g \left[\left(\frac{g}{f}\right)^\lambda-1\right],
~~~~ -\infty<\lambda<\infty.
\end{equation}
The family of divergences in (\ref{EQ:LDPD})
is similar to the density power divergence family given by
 \begin{equation}\label{EQ:dpd}
    {\rm DPD}_\alpha(g,f) = 
\int  \left[f^{1+\alpha} - \left(1 + 
\frac{1}{\alpha}\right)  f^\alpha g + \frac{1}{\alpha} g^{1+\alpha}\right] ~
{\rm for} ~\alpha \geq 0.
\end{equation}
The LDPD family may be recovered from the DPD by replacing the 
identity function with the logarithm function. 
In spite of the similarity between the forms of the 
DPD and the LDPD families, Jones et al.~(2001) had originally 
developed the latter as a special case of an estimation 
method proposed by Windham (1995). 
Following the connection shown
between PD and DPD by Patra et al.~(2013), 
Maji, Chakraborty and Basu (2014) have recently shown that the same sort of connection exists
between LDPD and the logarithmic power divergence (LPD) family.
To show this connection, the LDPD measure can be written as 
\begin{eqnarray}\label{EQ:LDPD_log}
{\rm LDPD}_\beta(g,f) &=& \displaystyle \log \int  f^{1+\beta} - \left(1 + \frac{1}{\beta}\right)  
\log \int \left(\frac{g}{f}\right) f^{1+\beta} \\ \nonumber
& & + \frac{1}{\beta} \log \int \left(\frac{g}{f}\right)^{1+\beta} f^{1+\beta}.
\end{eqnarray}
Replacing the $f^{1+\beta}$ term with $f$ 
in each of the three terms on the right hand side of (\ref{EQ:LDPD_log}) 
leads
to the density-based divergence $\frac{1}{\beta} \log \int \frac{g^{1+\beta}}{f^\beta}.$
After standardizing this further we express this divergence as 
 $\frac{1}{\beta(\beta+1)} \log \int \frac{g^{1+\beta}}{f^\beta}.$ This family of
divergences will be called the logarithmic power divergence family 
by Maji, Chakraborty and Basu (2014). Using a
different symbol for the tuning parameter, this
family has the form
\begin{equation}
\label{EQ:LPD}
{\rm LPD}_{\gamma}(g,f)=\frac{1}{\gamma(\gamma+1)} \log \int 
\frac{g^{1+\gamma}}{f^\gamma} , \gamma \in \mathbb{R}.
\end{equation}
The limiting forms as $\gamma \rightarrow 0$ and $\gamma \rightarrow -1$
generate, respectively, the likelihood disparity LD (or PD$_{0}$) and the 
Kullback-Leibler divergence KLD (or PD$_{-1}$). 
These divergences have the form 
\begin{equation}
 {\rm LD}(g,f)=  \int g \log \left(\frac{g}{f}\right),
\end{equation}
\begin{equation}
 {\rm KLD}(g,f)=  \int f \log \left(\frac{f}{g}\right).
\end{equation}
For any other value 
of $\gamma$ the LPD measure can be seen to be a function of the PD 
measure at the same value of $\gamma$. Specifically,
\begin{equation}
 \label{EQ:pd_lpd}
{\rm LPD}_\gamma(g,f)=\frac{1}{\gamma(\gamma+1)} \log 
\left[\gamma(\gamma+1){\rm PD}_\gamma(g,f)+1\right].
\end{equation}
Apart from being 
briefly considered by
Renyi as a measure of the amount of information
(Renyi, 1961), 
the LPD family
is a member 
of the ($h,\phi$) divergence family (eg. Pardo, $2006$),
where 
$h(x)=\frac{1}{r(r+1)} [r(r+1)\log x+1]$. 
\section{The Logarithmic Super Divergence}
\label{SEC:lsd_defn}
We now define the Logarithmic Super Divergence 
(or Logarithmic \\ $S$-Divergence)
introduced in Maji, Ghosh and Basu (2014)
and establish 
the asymptotic properties of the procedures resulting from it. 
The
Logarithmic $S$-Divergence (LSD) is defined (Maji, Ghosh and Basu, 2014) as
\begin{eqnarray}
\label{lsd_form}
		{\rm LSD}_{\beta, \gamma}(g,f) 
 &=&  \frac{1}{A} ~ \log \int f^{1+\beta}  -   \frac{1+\beta}{A B} ~ 
\log \int f^{B} g^{A} \nonumber \\ 
& & + \frac{1}{B} ~ \log \int g^{1+\beta},
~~~~ \beta > 0, ~-\infty < \gamma < \infty,
\end{eqnarray}
where
\begin{eqnarray}
    A &=& 1+\gamma (1-\beta) \nonumber \\ 
\mbox{           and               }
    B &=& \beta - \gamma (1-\beta). \nonumber
\end{eqnarray}
For $\beta=0~(A=1+\gamma,B=-\gamma),$ this family 
coincides with the 
logarithmic power divergence family of (\ref{EQ:LPD}) with parameter $\gamma$, 
while $\gamma=0$ gives the logarithmic density power divergence family
in (\ref{EQ:LDPD}) with parameter $\beta$. 
\paragraph{}
Replacing the logarithmic function with the identity function in (\ref{lsd_form})
generates the divergence
\begin{eqnarray}
		S_{\alpha, \lambda}(g,f)
 &=&  \frac{1}{A} ~ \int ~ f^{1+\alpha}  -   \frac{1+\alpha}{A B}
~ \int ~~ f^{B} g^{A}  + \frac{1}{B} ~ \int ~~ g^{1+\alpha},
\label{EQ:general_form}
\end{eqnarray}
with $\beta$ and $\gamma$ being replaced by $\alpha$ and $\lambda$
respectively. 
This divergence is known as the $S$-divergence and has been introduced by Ghosh et al.~(2013).
%
%
%
\subsection{Estimating Equation}
\label{SEC:lsd_est_eqn}
Consider a parametric class of model densities 
$\{ f_\theta : \theta \in \Theta \subseteq {\mathbb{R}}^p \}$
and 
suppose that our interest is in estimating 
$\theta$. Let $G$ denote the distribution function corresponding
to the true density $g$. The minimum 
LSD
functional
$T_{\beta,\gamma}(G)$ at $G$ is defined as 
\begin{equation}
 {\rm LSD}_{\beta, \gamma}\left(g,f_{T_{\beta,\gamma}(G)}\right)
= \min \limits_{\theta \in \Theta } {\rm LSD}_{\beta, \gamma}(g,f_\theta).
\end{equation}
It takes the value $\theta$ when the true density $g = f_\theta$ is in the model; when 
it does not, $\theta_{\beta,\gamma}^g = T_{\beta,\gamma}(G)$ represents 
the best fitting parameter, and $f_{\theta^g}$ is the model element closest 
to $g$ in terms of logarithmic super divergence. For simplicity of 
 notation, we suppress the subscript $\beta,\gamma$ in $\theta_{\beta,\gamma}^g$.
%
A simple differentiation gives us the estimating equation 
for $\theta$, which is 
\begin{equation}
\label{EQ:lsd_est_eqn}
 \frac{\int f_\theta^{1+\beta}  u_\theta}{\int f_\theta^{1+\beta}} = 
\frac{\int f_\theta^{B} g^A  u_\theta}{\int f_\theta^{B} g^A }.
\end{equation}
For $\beta=0~(A=1+\gamma,B=-\gamma)$, the equation thus becomes the same as the estimating 
equation of the logarithmic power divergence family with parameter $\gamma$. 
For $\gamma=0~(A=1, B=\beta)$, on the other hand, it is the estimating 
equation for the LDPD measure. 
Here $u_\theta(x)=\frac{\partial}{\partial \theta} \log f_\theta(x)$ 
is the likelihood score function. 

\section{Asymptotic Properties of the Minimum LSD Estimators in 
Discrete Models}
\label{SEC:lsd_properties}
Under the parametric set-up of Section {\ref{SEC:lsd_est_eqn}},
consider a discrete family of distributions. We will use 
the term ``density function'' generally for the sake of a 
unified notation, irrespective of whether the distribution is discrete
or continuous. Let $X_1,\ldots,X_n$ be a random sample from the 
true distribution having density function $g$. Representing 
the logarithmic $S$-divergence in terms of the parameter $\beta$ and $\gamma$
(as given in Section \ref{SEC:lsd_defn}), let 
$\hat{\theta}_{\beta,\gamma}$ be the estimator obtained by minimizing 
${\rm LSD}_{\beta,\gamma}(\hat{g},f_\theta)$ over $\theta \in \Theta$, where 
$\hat{g}$ is a suitable nonparametric density estimate of $g$; 
in the discrete case the vector of relative frequencies 
based on the sample data is the canonical choice for $\hat{g}$. 
For similar theoretical properties of the minimum $S$-divergence estimator, see Ghosh (2014).
\subsection{Asymptotic Properties} 
	Suppose $X_1,X_2,\ldots,X_n$ are $n$ independent 
and identically distributed observations from a discrete 
distribution $G$ modeled by 
$\mathcal{F} = \{ F_{\theta} : \theta \in \Theta \subseteq \mathbb{R}^p \}$ and 
let the distribution have support $\mathcal{X} = \{ 0,1,2,\ldots \}$ without loss
of generality. Denote the 
relative frequency at $x$ as obtained from the data by $r_n(x) = \frac{1}{n} 
\sum \limits_{i=1}^n \chi(X_i = x)$ where $\chi(A)$ is the 
indicator of $A$. 
The minimum LSD estimator 
is then obtained 
by minimizing the LSD 
measure between the data 
density $r_n$ and the model density $f_\theta$
with respect to $\theta$.

Minimizing
${\rm LSD}_{\beta,\gamma}(r_n,f_{\theta})$ 
with respect to $\theta$ is equivalent to
minimizing
$H_n(\theta)$ where
\begin{eqnarray}
    H_n(\theta) = \frac{1}{1+\beta} \left[ \frac{1}{A} \log \sum_x f_{\theta}^{1+\beta}(x) 
- \frac{1+\beta}{AB} \log \sum_x f_{\theta}^B(x) r_n^A(x) \right].
\end{eqnarray}
	Now,
\begin{eqnarray}
    \nabla H_n(\theta) &=& \frac{1}{A} 
\left[\frac{\sum_x f_{\theta}^{1+\beta}(x) u_\theta(x)}{\sum_x f_{\theta}^{1+\beta}(x)} -
\frac{\sum_x f_{\theta}^B(x) r_n^A(x) u_\theta(x)}{\sum_x f_{\theta}^B(x) r_n^A(x)}\right],
\end{eqnarray}
where $\nabla$ is 
the gradient with respect to $\theta$.
Equating 
the above
to zero, the
estimating equation becomes 
\begin{eqnarray}
    \sum_x K(\delta_n(x)) f_{\theta}^{1+\beta}(x) w_{\theta}(x) = 0,
\end{eqnarray}
where $\delta_n(x) = \frac{r_n(x)}{f_{\theta}(x)},~~
\delta_g(x) = \frac{g(x)}{f_{\theta}(x)},~~
w_{\theta}(x)=[B(\theta)u_\theta(x)-A(\theta)],~~
K(\delta)=\delta^A-1,~~ A(\theta)=\sum_x 
f_{\theta}^{1+\beta}(x)u_{\theta}(x) \\ \mbox{ and } 
B(\theta)=\sum_x 
f_{\theta}^{1+\beta}(x).$

	Define 
\begin{eqnarray}
J_g &=& E_g \left[ w_{\theta^g}(X)u_{\theta^g}^T(X) K'(\delta_g^g(X))
 f_{\theta^g}^{\beta}(X) \right] \nonumber \\
 & & -  \sum_x  K(\delta_g^g(x)) 
f_{\theta}^{1+\beta}(x) \nabla w_{\theta^g}(x) \nonumber \\
& & - (1+\beta) \sum_x  K(\delta_g^g(X)) 
f_{\theta}^{1+\beta}(x) w_{\theta^g}(x) u_{\theta^g}(x) . \\
    V_g &=& V_g \left[ K'(\delta_g^g(X)) 
f_{\theta^g}^\beta(X) w_{\theta^g}(X) \right].
\end{eqnarray}
We now present the necessary assumptions for our asymptotic results.  

\begin{enumerate}
    \item The model family ${F_\theta}$ is identifiable which 
means that different values of the parameter must generate 
different probability distributions of the observable variables.
	 \item The probability density function $f_\theta$ of 
the model distribution have common support so that the 
set $\mathcal{X} = \{ x : f_{\theta}(x) >0 \}$ is independent of $\theta$.
			Also the true distribution $g$ is 
compatible with the model family.
	 \item There exists an open subset $\omega \subset \Theta$
 for which the best fitting parameter $\theta^g$ is an interior 
point and for almost all $x$, the density $f_{\theta}(x)$ admits 
all third derivatives of the type $\nabla_{jkl} f_{\theta}(x)$   $\mbox{ for all } \theta \in \omega$. 
Here the subscripts $j,k,l$ of $\nabla$ represent the indicated 
partial derivatives. 
		\item The matrix $\frac{1+\beta}{A} J_g$ is positive definite.
		\item The quantities 
$\sum_x g^{1/2}(x) 
f_{\theta}^{\beta}(x) |u_{j\theta}(x)|$, 
$\sum_x g^{1/2}(x) 
f_{\theta}^{\beta}(x) |w_{j\theta}(x)|$, \\
$\sum_x g^{1/2}(x) 
f_{\theta}^{\beta}(x) |u_{j\theta}(x)| |u_{k\theta}(x)|$ and 

$\sum_x g^{1/2}(x) f_{\theta}^{\beta}(x) |u_{jk\theta}(x)|$  
are bounded $\mbox{ for all } j, k$ and $\mbox{ for all } \theta \in \omega$.

		\item For almost all $x$, there exists functions 
$M_{jkl}(x)$, $M_{jk,l}(x)$, $M_{j,k,l}(x)$,
possibly depending on $\beta$ that dominate, in absolute 
value, 
$f_{\theta}^{\beta}(x) u_{jkl\theta}(x)$, $f_{\theta}^{\beta}(x) 
u_{jk\theta}(x) u_{l\theta}(x)$ and $f_{\theta}^{\beta}(x) 
u_{j\theta}(x) u_{k\theta}(x) u_{l\theta}(x)$ respectively $ \mbox{ for all } j, k, l$ and 
that are uniformly bounded in expectation with respect to $g$ and $f_{\theta} 
~~ \mbox{for all } \theta \in \omega$, where 
the subscripts of $u$ denote the indicated partial 
derivatives of the score function.
		\item The function $\left( \frac{g(x)}{f_{\theta}(x)} \right)^{A-1}$ is 
uniformly bounded (by, say, $C$) $\mbox{ for all } \theta \in \omega$.

\end{enumerate}
	To prove the consistency and asymptotic normality of 
the minimum LSD estimator, we will assume, for the rest of the paper, that the 
seven conditions stated above are satisfied.
\paragraph{}
Based on the above assumptions we now 
start the proofs of the required results. 
The proofs are primarily along the 
lines of Lindsay (1994) and Basu et al.~(2011).
Define $ \eta_n(x) = \sqrt n (\sqrt
\delta_n-\sqrt\delta_g)^2$.
We then have the following lemmas. 
\begin{lemma}
\label{lem:lamma1}
  	 For any $k \in [0,2]$, we have
\begin{enumerate}
    \item $E[\eta_n^k(x)] \le n^{\frac{k}{2}}E[|\delta_n(X) - 
\delta_g(X)|]^k \le \left[ \frac{g(x)(1-g(x))}{f_{\theta}^2(x)}\right]^{\frac{k}{2}}.$ 
	 \item   $E[|\delta_n(X) - \delta_g(X)|] \le \frac{2g(x)(1-g(x))}{f_{\theta}(x)}.$
\end{enumerate}
\begin{proof}
For $a,b \ge 0$, we have the inequality $(\sqrt a - \sqrt b)^2 \le |a-b|$. So we get
\begin{eqnarray}
    E[\eta_n^k(x)] &=& n^{\frac{k}{2}} E[(\sqrt
\delta_n-\sqrt\delta_g)^2]^k \nonumber \\
		& \le & n^{\frac{k}{2}} E[|
\delta_n-\delta_g |]^k. \nonumber 
\end{eqnarray}

	For the next part see that, $ n r_n(x) \sim Bin( n, g(x) ) ~~ \mbox{ for all } x$. Now, 
for any  $k \in [0,2]$, we get by the Lyapunov's inequality that 

\begin{eqnarray}
    E[|\delta_n(X) - \delta_g(X)|]^k &\le & \left[ E(\delta_n(X) - \delta_g(X))^2 \right]^{\frac{k}{2}} \nonumber \\					&=& \frac{1}{f_{\theta}^k(x)} \left[ E(r_n(X) - g(X))^2 \right]^{\frac{k}{2}} \nonumber 
\\		&=& \frac{1}{f_{\theta}^k(x)} \left[ \frac{g(x)(1-g(x))}{n}\right]^{\frac{k}{2}}. \nonumber 
\\ \nonumber
\end{eqnarray}

	For the second part, note that
\begin{eqnarray}
    E[|\delta_n(X) - \delta_g(X)|] &=& \frac{1}{f_{\theta}^k(x)} \left[ E|r_n(X) - g(X)| \right]^{\frac{k}{2}} \nonumber 
\\		&\le& \frac{2g(x)(1-g(x))}{f_{\theta}(x)}. \nonumber 
\\ \nonumber
\end{eqnarray}
where the last inequality follows from the result about the mean-deviation of a Binomial random variable.
\end{proof}

\end{lemma}

\begin{lemma}
\label{lem:lamma2}
	$E[\eta_n^k(x)] \rightarrow 0$, as $n \rightarrow \infty$,  for $k \in [0,2).$
\begin{proof}
  This follows from Theorem 4.5.2 of Chung (1974) by noting that \\ 
$n^{1/4} (r_n^{1/2}(x) - g^{1/2}(x)) \rightarrow 0$ 
with probability one for each $x \in \mathcal{X}$ and by the Lemma \ref{lem:lamma1}(1), $sup_n E[\eta_n^k(x)] $ is bounded.
\end{proof}

\end{lemma}
%
	Let us now define, $a_n(x) = K(\delta_n(x)) - K(\delta_g(x))$ and 
$b_n(x) = (\delta_n(x)-\delta_g(x))K'(\delta_g(x))$.

 We will need the limiting distributions of 
$$S_{1n} = \sqrt n \sum_x a_n(x)f_{\theta}^{1+\beta}(x)w_{\theta}(x)$$
and $$S_{2n} = \sqrt n \sum_x 
b_n(x)f_{\theta}^{1+\beta}(x)w_{\theta}(x).$$

 Define $\tau_n(x) = \sqrt n |a_n(x) - b_n(x)|$. 

\begin{lemma}
\label{lem:lamma3}
	Suppose Assumption 5 holds. Then $E|S_{1n}-S_{2n}| \rightarrow 0 $ as $n \rightarrow \infty$,
and hence $S_{1n}-S_{2n} \xrightarrow{P} 0 $ as $n \rightarrow \infty.$
\begin{proof}
 	By Lemma 2.15 of Basu et ~al.~(2011) [or, Lindsay(1994), Lemma 25], there 
exists some positive constant $\beta$ such that 
$$
\tau_n(x) \le \beta \sqrt n (\sqrt\delta_n-\sqrt\delta_g)^2 = \beta \eta_n(x).
$$ 
Also, by Lemma \ref{lem:lamma1}, $ E[\tau_n(x)] \le \beta \frac{g^{1/2}(x)}{f_{\theta}(x)}$.\\
And by Lemma \ref{lem:lamma2}, $E[\tau_n(x)] = \beta E[\eta_n(x)] \rightarrow 0$ as $n \rightarrow \infty$. Thus we get, 
\begin{eqnarray}
    E|S_{1n}-S_{2n}| &\le& \sum_x E[\tau_n(x)] f_{\theta}^{1+\beta}(x)|w_{\theta}(x)| \nonumber \\
                     &\le&  \beta \sum_x g^{1/2}(x) f_{\theta}^{\beta}(x)|w_{\theta}(x)| \nonumber \\
							& < &   
\infty ~~~~~~ \mbox{ (by Assumption $5$)}. \nonumber \\ \nonumber
\end{eqnarray}
So, by the Dominated Convergence Theorem (DCT),  $E|S_{1n}-S_{2n}| \rightarrow 0 $ as $n \rightarrow \infty$.
Hence, by Markov's inequality, $S_{1n}-S_{2n} \xrightarrow{P} 0 $ as $n \rightarrow \infty$.
\end{proof}

\end{lemma}
%
%
\begin{lemma}
\label{lem:lamma4}
	Suppose $V_g =  V_g \left[ K'(\delta_g(X)) f_{\theta}^\beta(X) w_{\theta}(X) \right] $ is finite. Then 
$$
S_{1n} \rightarrow N(0, V_g).
$$
\begin{proof}
 	By Lemma \ref{lem:lamma3}, the asymptotic distribution of 
$S_{1n}$ and $S_{2n}$ are the same. Now, we have
\begin{eqnarray}
    S_{2n} &=& \sqrt n \sum_x (\delta_n(x) - \delta_g(x)) 
K'(\delta_g(x)) f_{\theta}^{1+\beta}(x) w_{\theta}(x) \nonumber \\
				&=&  \sqrt n \sum_x (r_n(x) - 
g(x)) K'(\delta_g(x)) f_{\theta}^{\beta}(x) w_{\theta}(x) \nonumber \\
			&=&  \sqrt n \left( \frac{1}{n} 
\sum_{i=1}^n \left[ K'(\delta_g(X_i)) f_{\theta}^{\beta}(X_i) w_{\theta}(X_i) - 
E_g\{ K'(\delta_g(X))f_{\theta}^{\beta}(X)w_{\theta}(X)\} \right]\right) \nonumber \\
			&\rightarrow & N( 0, V_g) ~~~~~ \mbox{(by the
Central Limit Theorem)}. \nonumber \\ \nonumber 
\end{eqnarray}
\end{proof}

\end{lemma}
%
	We will now consider the final theorem of this section 
about the consistency and asymptotic normality of the minimum LSD estimator.
\begin{theorem}
\label{th:est_th}
	Under Assumptions $1-7$, there exists a consistent 
sequence $\theta_n$ of roots to the minimum LSD estimating equation (\ref{EQ:lsd_est_eqn}).
	Also, the asymptotic distribution of 
$\sqrt n (\theta_n - \theta_g)$ is $p-$dimensional normal 
with mean $0$ and variance $J_g^{-1}V_g J_g^{-1}$.
\begin{proof}
Because of the lengthy and somewhat messy calculations, the 
proof of consistency has been put in Appendix.
\paragraph{}
{\bf Proof of the asymptotic Normality :}		For the Asymptotic normality, we expand 
$$\sum_x K(\delta_n(x)) f_{\theta}^{1+\beta}(x)w_{\theta}(x)$$
 in Taylor series about $\theta = \theta^g$ to get
\begin{eqnarray}
\label{EQ:kfw_exapnd}
    & & \sum_x K(\delta_n(x)) f_{\theta}^{1+\beta}(x)w_{\theta}(x)  \nonumber \\
& = & \sum_x K(\delta_n^g(x)) f_{\theta^g}^{1+\beta}(x)w_{\theta^g}(x) 
+ \sum_k (\theta_k - \theta_k^g) \nabla_k \left( \sum_x K(\delta_n(x)) 
f_{\theta}^{1+\beta}(x)w_{\theta}(x) \right)|_{\theta = \theta^g} \nonumber \\
 & & ~~~~~~~~~ + \frac{1}{2} \sum_{k,l} (\theta_k - \theta_k^g)
(\theta_l - \theta_l^g)\nabla_{kl} \left(\sum_x K(\delta_n(x)) 
f_{\theta}^{1+\beta}(x)w_{\theta}(x)\right)|_{\theta = \theta '} \\ \nonumber
\end{eqnarray}
where, $\theta '$ lies in between $\theta$ and $\theta^g$.
	Now, let $\theta_n$ be the solution of the minimum 
LSD estimating equation, which can be assumed to be consistent. 
Replace $\theta$ by $\theta_n$ in above (\ref{EQ:kfw_exapnd}) 
so that the LHS of the equation becomes zero and hence we get
\begin{eqnarray}
\label{EQ:kfw_times}
  &  & - \sqrt n \sum_x K(\delta_n^g(x)) f_{\theta^g}^{1+\beta}(x)
w_{\theta^g}(x)  =  \sqrt n \sum_k (\theta_{nk} - 
\theta_k^g) \times \nonumber \\
 & & ~~~~~~~~~~~~ \{ \nabla_k \left( \sum_x K(\delta_n(x)) 
f_{\theta}^{1+\beta}(x)w_{\theta}(x) \right)|_{\theta = \theta^g} 
+ \nonumber \\
& & ~~~~~~~~~~~~~~~~~~~~~~~\frac{1}{2} \sum_{l} 
(\theta_{nl} - \theta_l^g)\nabla_{kl} \left(\sum_x K(\delta_n(x)) 
f_{\theta}^{1+\beta}(x)w_{\theta}(x)\right)|_{\theta = \theta '}  
\}. \nonumber \\
\end{eqnarray}
	Note that, the first term within the bracketed quantity 
in the RHS of above (\ref{EQ:kfw_times}) converges to $J_g$ with probability 
tending to one, while the second bracketed term is an 
$o_p(1)$ term (as proved in the proof of consistency part). Also, 
by using the Lemma \ref{lem:lamma4}, we get that	
\begin{eqnarray}
   & & \sqrt n \sum_x K(\delta_n^g(x)) 
f_{\theta^g}^{1+\beta}(x)w_{\theta^g}(x) \nonumber \\
 ~~~~~~~~~~~ & = & \sqrt n \sum_x  [K(\delta_n^g(x)) - 
K(\delta_g^g(x))]  f_{\theta^g}^{1+\beta}(x)w_{\theta^g}(x)  (\mbox{Using estimating equation)} \nonumber \\
 ~~~~~~~~~~~ & = & S_{1n}|_{\theta=\theta^g} 
\rightarrow N_p(0, V_g). \nonumber \\ 
\end{eqnarray}
Therefore, by Lehmann(1983, Lemma 4.1), 
$\sqrt n(\theta_n - \theta^g)$ has asymptotic 
distribution as $ N_p(0,J_g^{-1} V_g J_g^{-1})$. 
\end{proof}
\end{theorem}
A very interesting observation that follows from the 
asymptotic distribution just established is that the 
asymptotic distribution of the estimator is independent 
of the parameter $\gamma$. Yet the behavior of the 
estimator varies widely with $\gamma$ under the presence of
outliers. Here we briefly report the findings 
reported by Maji, Ghosh and Basu (2014) in this connection, 
which is at least partially indicated by the results 
of the current and the subsequent section. Maji, Ghosh and Basu (2014) 
have observed that the first order influence function of the minimum 
LSD estimator is independent of $\gamma$, predicting that the 
robustness properties of the minimum LSD estimators are 
similar for each value of $\gamma$. This is immediately 
contradicted by the other results of Maji, Ghosh and Basu (2014)
as well as the next section of the current article. 
Further exploration by Maji, Ghosh and Basu (2014) indicate 
that the second order influence function 
gives a much more accurate picture of the 
robustness of these estimators. This clearly indicates the 
limitation of the first order influence function in quantifying 
the robustness of the estimators in this case. In fact 
the second order influence analysis (see Maji, Ghosh and Basu, 2014 for details)
shows that the limitation of the first order influence function 
can go both ways -- it can fail to indicate the stability of a robust 
estimator, and can also describe a highly unstable estimator as 
a robust one.
\section{Examples}
\label{SEC:est_example}
\subsection{An Example with a Geometric Model}
The data set given in Basu et al.~(2011),
Table $2.4$, represent the cases of peritonitis
for $390$ kidney patients.
Following Basu et al.~(2011) we 
use a 
geometric model with 
%
parameter $\theta$ (success probability) 
as our target distribution. 
A quick look at the data reveals that a geometric model
with $\theta=0.5$ may fit the data well. 
We fit a geometric model with parameter $\theta$ using the 
LSD measure for various combinations of $\beta$ and $\gamma$. 
We can treat the two observations at 10 
and 12 as mild to moderate outliers. We have 
evaluated the minimum
LSD estimator 
in this case for the full data as well as the outlier 
deleted data.
The estimates are 
presented 
shown in
Tables \ref{tab:kidney_w_out} and \ref{tab:kidney_wo_out} respectively. 
The estimates highlight an interesting point; 
for 
$\gamma<0$ or for larger values of $\beta$ with $\gamma\geq0$ the 
parameter estimates are close for full data and outlier deleted data.
However for 
$\gamma>0$ and 
$\beta$ small, the presence or absence 
of the outliers do not lead to a substantially
larger difference. 
This gives a clear indication about which 
combinations of the $(\beta,\gamma)$ values keep 
the estimators stable and which are the ones that are easily 
affected. 

 \begin{table}

	\centering 
	\caption{The estimates of the parameter 
of the geometric model for 
different values of $\gamma$ and $\beta$ for the 
Peritonitis Incidence Data with outlier
}
		\begin{tabular}{c c c c c c c } \hline
$\gamma$ \textdownarrow $\beta$ $\rightarrow$ 
& 0 & 0.2  & 0.4  &  0.6 & 0.8 & 1 \\ 
\hline
$-1$ & -- & 0.518 & 0.511 & 0.511 & 0.513 & 0.515 \\ 
$-0.7$ & 0.519 & 0.51 & 0.509 & 0.51 & 0.512 & 0.515 \\ 
$-0.5$ & 0.510 & 0.508 & 0.508 & 0.51 & 0.512 & 0.515 \\ 
$-0.3$ & 0.504 & 0.505 & 0.507 & 0.51 & 0.512 & 0.515 \\ 
$-0.1$ & 0.499 & 0.503 & 0.506 & 0.509 & 0.512 & 0.515 \\ 
0 & 0.496 & 0.502 & 0.506 & 0.509 & 0.512 & 0.515 \\ 
0.5 & 0.48 & 0.496 & 0.504 & 0.508 & 0.512 & 0.515 \\ 
1 & 0.461 & 0.486 & 0.501 & 0.507 & 0.512 & 0.515 \\ 
1.3 & 0.45 & 0.479 & 0.499 & 0.507 & 0.511 & 0.515 \\ 
1.7 & 0.438 & 0.469 & 0.495 & 0.506 & 0.511 & 0.515 \\ 
2 & 0.43 & 0.461 & 0.491 & 0.505 & 0.511 & 0.515 \\ 
\hline
		\end{tabular}
\label{tab:kidney_w_out}
\end{table}
 \begin{table}

	\centering 
	\caption{The estimates of the parameter 
of the geometric model for 
different values of $\gamma$ and $\beta$ for the 
Peritonitis Incidence Data without the two outliers 
}
		\begin{tabular}{c c c c c c c } \hline
$\gamma$ \textdownarrow $\beta$ $\rightarrow$ 
& 0 & 0.2  & 0.4  &  0.6 & 0.8 & 1 \\ 
\hline
$-1$ & -- & 0.521 & 0.512 & 0.511 & 0.513 & 0.515 \\ 
$-0.7$ & 0.526 & 0.513 & 0.51 & 0.511 & 0.513 & 0.515 \\ 
$-0.5$ & 0.518 & 0.511 & 0.509 & 0.51 & 0.512 & 0.515 \\ 
$-0.3$ & 0.513 & 0.509 & 0.508 & 0.51 & 0.512 & 0.515 \\ 
$-0.1$ & 0.510 & 0.508 & 0.508 & 0.51 & 0.512 & 0.515 \\ 
0 & 0.509 & 0.507 & 0.508 & 0.51 & 0.512 & 0.515 \\ 
0.5 & 0.505 & 0.505 & 0.506 & 0.509 & 0.512 & 0.515 \\ 
1 & 0.501 & 0.503 & 0.505 & 0.508 & 0.512 & 0.515 \\ 
1.3 & 0.5 & 0.501 & 0.504 & 0.508 & 0.512 & 0.515 \\ 
1.7 & 0.498 & 0.5 & 0.504 & 0.508 & 0.511 & 0.515 \\ 
2 & 0.496 & 0.499 & 0.503 & 0.507 & 0.511 & 0.515 \\ 
\hline
		\end{tabular}
\label{tab:kidney_wo_out}
\end{table}
\subsection{An Example with a Poisson Model}
This example 
gives us the observed frequencies and 
corresponding estimated frequencies (Table 
\ref{droso_table1}) for several 
minimum LSD 
estimators 
under the Poisson model
for a sex linked recessive lethal test in drosophila (fruit flies)
exposed to a certain chemical. For each of several male flies
one samples about $100$ daughter flies, and 
then determines the frequency of the number of 
daughter flies having a recessive lethal mutation 
in its $X$-chromosome. The data represent a frequency 
of frequencies; refer to 
Woodruff et al.~(1984) for details. There is a possible case
of outliers corresponding to the observations at $x=3,4$. 
\paragraph{}
Table \ref{droso_table1} provides the estimators and the predicted 
frequencies for a small number of $\gamma,\beta$ combinations,
together with the fits of the maximum likelihood estimator
(denoted by ML) and the outlier deleted maximum likelihood
estimator obtained by removing the two outliers (denoted
by ML$+$D). Clearly the estimators (and the estimated frequencies)
are substantially different for the ML 
and ML$+$D cases, demonstrating that the maximum 
likelihood estimator is significantly affected by the presence 
of these outliers. Also apart from the ML, the ($\gamma=1,\beta=0.1$)
combination leads to highly unstable estimators.  
While Table \ref{droso_table1} provides a small number of ($\gamma,\beta$) 
combinations, a large selection is presented in Table \ref{droso_estimate}, 
where the salient features may be described as follows: 
\begin{enumerate}[(a)]
 \item All the estimators corresponding to large negative values of 
$\gamma$ and/or values of $\beta$ close to $1$ generate outlier 
resistant methods,
\item Estimators corresponding to large positive values of $\gamma$
are relatively poor in terms of robustness, especially for small $\beta$
\item All estimators for $\beta=1$ are identical,
irrespective of the value of $\gamma$. 
\end{enumerate}
Once again this example shows that large positive values of $\gamma$ 
with $\beta$ close to zero are the more unstable distances 
within the LSD class. 
\begin{table}
\caption{ 
Fits of the Poisson model to the Drosophila Data using several estimation methods
}
\vspace{2 mm}
\centering 
  \begin{tabular}{l c c c c c c c} 
\hline 
& \multicolumn{7}{c}{Recessive lethal count}\\\cline{2-8}
& 0 & 1 & 2 & 3 & 4 & $\geq 5$ 
& $\hat{\theta}$
\\
Observed & 23 & 3 & 0 & 1& 1& 0\\
ML & 19.59 & 7.00 & 1.25 & 0.15 & 0.01 & -- & 0.3571\\
ML+D & 24.95 & 2.88 & 0.17 & 0.01 & -- & -- & 0.1154\\
{\rm LSD}$_{\gamma = 1~~\&~~\beta = 0.1}$ & 14.90 & 9.40 & 2.97 & 0.62 & 0.10 & 0.01 & 0.6311\\
{\rm LSD}$_{\gamma = -1~~\&~~\beta = 0.1}$ & 25.95 & 1.98 & 0.07 & -- & -- & -- & 0.0762\\
{\rm LSD}$_{\gamma = 1~~\&~~\beta = 1}$ & 24.59 & 3.19 & 0.21 & 0.01 & -- & -- & 0.1297\\
{\rm LSD}$_{\gamma =-1~~\&~~\beta = 1}$ & 24.59 & 3.19 & 0.21 & 0.01 & -- & -- & 0.1297\\
\hline
\end{tabular}
\label{droso_table1}
 \end{table}
 \begin{table}
	\centering 
	\caption{Estimates of the parameter for drosophila data in Table {\ref{droso_table1}}}
		\begin{tabular}{c c c c c c c c c c c c} \hline
$\gamma$ \textdownarrow $\beta$ $\rightarrow$ 
& 0 & 0.2  & 0.4  &  0.6 & 0.8 & 0.9 & 1 \\ 
\hline
$-0.8$ & 0.088 & 0.113 & 0.123 & 0.127 & 0.129 & 0.129 & 0.13 \\ 
$-0.7$ & 0.101 & 0.117 & 0.124 & 0.127 & 0.129 & 0.129 & 0.13 \\ 
$-0.6$ & 0.112 & 0.121 & 0.126 & 0.128 & 0.129 & 0.129 & 0.13 \\ 
$-0.5$ & 0.125 & 0.126 & 0.127 & 0.128 & 0.129 & 0.13 & 0.13 \\ 
$-0.4$ & 0.145 & 0.132 & 0.129 & 0.129 & 0.129 & 0.13 & 0.13 \\ 
$-0.3$ & 0.177 & 0.139 & 0.131 & 0.129 & 0.129 & 0.13 & 0.13 \\ 
$-0.2$ & 0.227 & 0.151 & 0.134 & 0.13 & 0.13 & 0.13 & 0.13 \\ 
$-0.1$ & 0.291 & 0.169 & 0.137 & 0.131 & 0.13 & 0.13 & 0.13 \\ 
0 & 0.357 & 0.194 & 0.142 & 0.132 & 0.13 & 0.13 & 0.13 \\ 
0.1 & 0.417 & 0.228 & 0.148 & 0.133 & 0.13 & 0.13 & 0.13 \\ 
0.2 & 0.47 & 0.269 & 0.157 & 0.134 & 0.13 & 0.13 & 0.13 \\ 
0.3 & 0.514 & 0.311 & 0.169 & 0.135 & 0.13 & 0.13 & 0.13 \\ 
0.4 & 0.553 & 0.353 & 0.184 & 0.137 & 0.131 & 0.13 & 0.13 \\ 
0.5 & 0.586 & 0.393 & 0.204 & 0.139 & 0.131 & 0.13 & 0.13 \\ 
0.6 & 0.615 & 0.43 & 0.226 & 0.142 & 0.131 & 0.13 & 0.13 \\ 
0.7 & 0.641 & 0.463 & 0.252 & 0.145 & 0.131 & 0.13 & 0.13 \\ 
0.8 & 0.663 & 0.494 & 0.278 & 0.149 & 0.131 & 0.13 & 0.13 \\ 
\hline
		\end{tabular}
\label{droso_estimate}	 
\end{table}

\section{Testing Parametric Hypothesis using the LSD Measures}
\label{SEC:testing}

Now we focus our attention on hypothesis testing, 
the other important \\ paradigm of statistical inference. 
\subsection{One Sample problem}
	We consider a parametric family of densities 
$\mathcal{F} = \{f_\theta : \theta \in \Theta \subseteq \mathbb{R}^p\}$ 
as above. Suppose we are given a random sample $X_1, \ldots, X_n$ of 
size $n$ from the population. Based on this sample, we want to test the hypothesis 
$$
 H_0 : \theta = \theta_0 ~~~ \mbox{against} ~~~~ H_1 : \theta \ne \theta_0.
$$
When the model is correctly 
specified and the null hypothesis is correct, $f_{\theta_0}$ is the data 
generating density. 
We consider the 
test statistics based on the LSD with parameter 
$\beta$ and $\gamma$ 
as follows:
\begin{equation}
 W_{\beta,\gamma}(\hat{\theta}_{\beta,\gamma},\theta_0)
=2n~{\rm LSD}_{\beta,\gamma}(f_{\hat{\theta}_{\beta,\gamma}},f_{\theta_0}),
\end{equation}
where ${\rm LSD}_{\beta,\gamma}(f_{\hat{\theta}_{\beta,\gamma}},f_{\theta_0})$ has the form 
given in (\ref{lsd_form}). 
\begin{theorem}
 The 
asymptotic distribution of the test statistic 	$W_{\beta,\gamma}
(f_{{\hat{\theta}}_{\beta,\gamma}},f_{{\theta_0}})$, under the null hypothesis 
$H_0 : \theta = \theta_0$, coincides with the distribution of 
$$
\sum_{i=1}^r ~  \zeta_i^{\beta}(\theta_0)Z_i^2
$$
where $Z_1, \ldots,Z_r$ are independent standard normal variables, 
$\zeta_1^{\beta}(\theta_0), \ldots, \\ \zeta_r^{\beta}
(\theta_0)$  are the nonzero eigenvalues of $A_{\beta}(\theta_0) 
J_{\beta}^{-1}(\theta_0) K_{\beta}(\theta_0) J_{\beta}^{-1}(\theta_0)$, 
with $J_{\beta}(\cdot)$ and $K_{\beta}(\cdot)$ as defined in Theorem
\ref{th:est_th}
and the matrix $A_{\beta}(\theta_0)$ is defined as
$$
A_{\beta}(\theta_0) = \nabla [  \nabla {\rm LSD}_{\beta,\gamma}(f_\theta, f_{\theta_0}) ]
|_{\theta = \theta_0} 
$$
 and
$$
r = rank\left( J_{\beta}^{-1}(\theta_0) 
K_{\beta}(\theta_0) J_{\beta}^{-1}(\theta_0)
A_{\beta}(\theta_0) J_{\beta}^{-1}(\theta_0) 
K_{\beta}(\theta_0) J_{\beta}^{-1}(\theta_0)\right).
$$
\begin{proof}
  We consider the second order Taylor series expansion 
of ${\rm LSD}_{\beta,\gamma}(f_\theta, f_{\theta_0})$ around 
$\theta = \theta_0$ at $\theta = \hat{\theta}_{\beta}$ as,
\begin{eqnarray}
\label{EQ:lsd_expansion}
    {\rm LSD}_{\beta,\gamma}\left(f_{\hat{\theta}_{\beta}},f_{\theta_0}\right) 
&=& {\rm LSD}_{\beta,\gamma}(f_{\theta_0}, f_{\theta_0}) + \sum_{i=1}^p
 \nabla_i {\rm LSD}_{\beta,\gamma}(f_\theta, 
f_{\theta_0})|_{\theta = \theta_0} (\hat{\theta}_{\beta}^i - \theta_0^i) \nonumber 
\\
& & + \frac{1}{2} \sum_{i,j} \nabla_{ij} {\rm LSD}_{\beta,\gamma}
(f_\theta, f_{\theta_0})|_{\theta = \theta_0} 
(\hat{\theta}_{\beta}^i - \theta_0^i) (\hat{\theta}_{\beta}^j - 
\theta_0^j) \nonumber \\ 
& & + o(||\hat{\theta}_{\beta} - \theta_0||^2), 
\end{eqnarray}
where 
$\nabla_i$ and $\nabla_{ij}$ represent the indicated partial 
derivatives with respect to the components of $\theta$. 
Now 
we have 
$$ {\rm LSD}_{\beta,\gamma}(f_{\theta_0}, f_{\theta_0}) = 0 $$ 
and $$\nabla_i {\rm LSD}_{\beta,\gamma}(f_\theta, 
f_{\theta_0})|_{\theta = \theta_0} =0. $$ 
Note that the above second order partial derivative of 
${\rm LSD}_{\beta,\gamma}(f_\theta, f_{\theta_0})$ at 
$\theta = \theta_0 $ is independent of $\gamma$ and 
so we will denote that as function of $\beta$ only.
We will denote the second order partial derivatives of
${\rm LSD}_{\beta,\gamma}(f_\theta, f_{\theta_0})$  in
(\ref{EQ:lsd_expansion})
by $a_{ij}^{\beta} (\theta_0)$. 
Also denote
$A_{\beta}(\theta_0) = \left( a_{ij}^{\beta} (\theta_0) 
\right)_{i,j=1,\ldots,p}$.
Now from the above Taylor series expansion it is clear 
that the random variables 
\begin{eqnarray}
     W_{\beta,\gamma}({\hat{\theta}_{\beta,\gamma}},{\theta_0}) = 
2 n {\rm LSD}_{\beta,\gamma}(f_{\hat{\theta}_{\beta,\gamma}},f_{\theta_0}) \nonumber \\
~~~~ \mbox{and} ~~~~ \sqrt{n} (\hat{\theta}_{\beta,\gamma} - \theta_0)^T 
A_{\beta} (\theta_0) \sqrt{n} (\hat{\theta}_{\beta,\gamma} - \theta_0) \nonumber
\end{eqnarray}
%
%
%
have the same asymptotic distribution. Now we know from 
the previous section
that the 
asymptotic distribution of $\sqrt{n} (\hat{\theta}_{\beta,\gamma}-\theta_0)$ 
is normal with mean zero and variance $J_{\beta}^{-1}(\theta_0) 
K_{\beta}(\theta_0) J_{\beta}^{-1}(\theta_0)$.
   Further we know that for $X \sim N_q(0,\Sigma)$, and a $q-$dimensional 
real symmetric matrix $A$, the distribution of the quadratic form 
$X^T A X$ is the same as that of
$\sum \limits_{i=1}^r \zeta_i^{\beta} Z_i^2$, where $Z_1, \ldots ,Z_r$ are independent 
standard normal variables, $r = rank(\Sigma A \Sigma)$, $r \ge 1$ and 
$\zeta_1^{\beta}, \ldots , \zeta_r^{\beta}$ are the nonzero eigenvalues of $A \Sigma$ 
(Dik and Gunst, 1985, Corollary 2.1). Applying this result with 
$ X =  \sqrt{n} (\hat{\theta}_{\beta,\gamma} - \theta_0)$ 
the theorem is established. It
is evident that the asymptotic distribution of the statistic depends on $\beta$ only
and is independent of $\gamma$.
\end{proof}
\end{theorem}

\begin{theorem}
An approximation to the power 
function of the test statistic $W_{\beta,\gamma}( {\hat{\theta}_{\beta,\gamma}} 
, {\theta_0})$ for testing $H_0 : \theta =\theta_0$ against 
$H_1: \theta \ne \theta_0$ at the significance level $\alpha$ is given by
\begin{eqnarray}
    \pi_{n,\alpha}^{\beta,\gamma} (\theta^*) = 
1 - \Phi \left( \frac{\sqrt{n}}{\sigma_{\beta,\gamma 
}(\theta^*)} \left( \frac{t_\alpha^{\beta,\gamma}}{2 n} 
- {\rm LSD}_{\beta,\gamma}(f_{\theta^*}, f_{\theta_0}) \right) \right), ~~~ \theta^* \ne \theta_0
\end{eqnarray}
where $t_\alpha^{\beta,\gamma}$ is the $(1-\alpha)^{th}$ 
quantile of the asymptotic distribution of \\ $W_{\beta,\gamma}
( {\hat{\theta}_{\beta,\gamma}} , {\theta_0})$, and $\sigma_{\beta,\gamma}(\theta^*)$ is defined as
\begin{eqnarray}
\label{EQ:sigma_defn}
    \sigma_{\beta,\gamma}^2(\theta) = 
M_{\beta,\gamma}(\theta)^T J_\beta^{-1}(\theta) 
K_\beta(\theta) J_\beta^{-1}(\theta) M_{\beta,\gamma}(\theta)
\end{eqnarray}
with 
$$
M_{\beta,\gamma}(\theta) = \nabla {\rm LSD}_{\beta,\gamma}
( f_\theta , f_{\theta_0}). 
$$
\begin{proof}
%
Fix some $\theta^* \ne \theta_0$. Consider 
the first order Taylor series expansion of 
${\rm LSD}_{\beta,\gamma}( f_{\hat{\theta}_{\beta,\gamma}} , 
f_{\theta_0})$ under $f_{\theta^*}$ as
$$
{\rm LSD}_{\beta,\gamma}( f_{\hat{\theta}_{\beta,\gamma}} , 
f_{\theta_0}) = {\rm LSD}_{\beta,\gamma}( f_{\theta^*} , 
f_{\theta_0}) + M_{\beta,\gamma}(\theta^*)^T 
(\hat{\theta}_{\beta,\gamma} - \theta^*) + o(||\hat{\theta}_{\beta,\gamma} - \theta^*||)
$$
where $M_{\beta,\gamma}$ is as defined in the theorem. Now 
we know that, under $\theta^*$, 
$$
\sqrt{n} (\hat{\theta}_{\beta,\gamma} - \theta^*) \rightarrow 
N(0, J_\beta^{-1}(\theta^*) K_\beta(\theta^*) J_\beta^{-1}(\theta^*) )
 ~~~as ~ n \rightarrow \infty
$$
and $ \sqrt{n} o(||\hat{\theta}_{\beta,\gamma} - \theta^*||) = o_p(1)$. Thus 
we get that the random variables 
$$
\sqrt{n} \left[ {\rm LSD}_{\beta,\gamma}( f_{\hat{\theta}_{\beta,\gamma}} , 
f_{\theta_0}) - {\rm LSD}_{\beta,\gamma}( f_{\theta^*} , f_{\theta_0}) 
\right] ~~~~~\mbox{ and } ~~~~~ M_{\beta,\gamma}(\theta^*)^T 
\sqrt{n}(\hat{\theta}_{\beta,\gamma} - \theta^*) 
$$
have the same asymtotic distribution . Therefore, we have 
$$
\sqrt{n} \left[ {\rm LSD}_{\beta,\gamma}( f_{\hat{\theta}_{\beta,\gamma}} , 
f_{\theta_0}) - {\rm LSD}_{\beta,\gamma}( f_{\theta^*} , f_{\theta_0}) 
\right] \rightarrow N( 0, \sigma_{\beta,\gamma}(\theta^*) )
$$
where $\sigma_{\beta,\gamma}(\theta^*)$ is as given in 
(\ref{EQ:sigma_defn}) above. Hence the desired approximation to the power function 
follows from the above asymptotic distribution. 
\end{proof}

\end{theorem}

\subsection{Two-Sample Problem}
	Again consider a parametric family of densities 
$\{f_\theta : \theta \in \Theta \subseteq \mathbb{R}^p\}$ as above in 
one sample problem, but here we are given two random 
samples $X_1, \ldots, X_n$ of size $n$ and $Y_1, \ldots, Y_m$ of 
size $m$ from two 
distributions within the model
having parameters $\theta_1$ and 
$\theta_2$ respectively and based on these two samples, we want to 
test for the homogeneity of the two samples, i.e. to test the hypothesis 
$$
 H_0 : \theta_1 = \theta_2 ~~~ \mbox{against} ~~~~ H_1 : \theta_1 \ne \theta_2.
$$
 We will consider the estimator $^{(1)}\hat{\theta}_{\beta,\gamma}$  and 
$^{(2)}\hat{\theta}_{\beta,\gamma}$ of $\theta_1$ and $\theta_2$ respectively, 
obtained by minimizing the LSD having parameters ${\beta,\gamma}$ and,
as before, 
will 
consider the test statistic based on the LSD with parameter $\beta$ 
and $\gamma$ as 
given by 
\begin{eqnarray}
    S_{\beta,\gamma}\left( ^{(1)}\hat{\theta}_{\beta,\gamma} , 
^{(2)}\hat{\theta}_{\beta,\gamma}\right) 
= \frac{2nm}{n+m} ~  {\rm LSD}_{\beta,\gamma}\left( f_{^{(1)}\hat{\theta}_{\beta,\gamma}} , 
f_{^{(2)}\hat{\theta}_{\beta,\gamma}}\right).
\end{eqnarray}
Now, first let us 
consider the asymptotic distribution of the test statistic $S_{\beta,\gamma}
\left( ^{(1)}\hat{\theta}_{\beta,\gamma}, ^{(2)}\hat{\theta}_{\beta,\gamma}\right)$ under $H_0$ 
in the following theorem. 
Assume that $\frac{m}{m+n}\rightarrow \omega, (0<\omega<1), \mbox{ as } m\rightarrow \infty
\mbox{ and } n\rightarrow\infty$. 
\begin{theorem}
\label{th_two}
The asymptotic distribution of the test 
statistic \\	$S_{\beta,\gamma}\left( ^{(1)}\hat{\theta}_{\beta,\gamma} , 
^{(2)}\hat{\theta}_{\beta,\gamma}\right)$, under the null hypothesis $H_0 : \theta_1 = \theta_2$, 
coincides with the distribution of 
$$
\sum_{i=1}^r ~  \zeta_i^{\beta}(\theta_1)Z_i^2
$$
where $Z_1, \ldots,Z_r$ are independent standard normal variables, 
$\zeta_1^{\beta}(\theta_1), \ldots, \\ \zeta_r^{\beta}(\theta_1)$ 
are the nonzero eigenvalues of $A_{\beta}(\theta_1) J_{\beta}^{-1}(\theta_1) 
K_{\beta}(\theta_1) J_{\beta}^{-1}(\theta_1)$, with \\
$J_{\beta}(\cdot)$, $K_{\beta}(\cdot)$ 
and $A_{\beta}(\cdot)$ as defined in the previous section and
$$
r = rank\left( J_{\beta}^{-1}(\theta_1) 
K_{\beta}(\theta_1) J_{\beta}^{-1}(\theta_1)
A_{\beta}(\theta_1) J_{\beta}^{-1}(\theta_1) 
K_{\beta}(\theta_1) J_{\beta}^{-1}(\theta_1)\right).
$$
\begin{proof} 
 We have 
$$
\sqrt{n} \left( ^{(1)}\hat{\theta}_{\beta,\gamma} - 
\theta_1 \right) \rightarrow N( 0, J_{\beta}^{-1}
(\theta_1) K_{\beta}(\theta_1) J_{\beta}^{-1}(\theta_1) )
$$
and
$$
\sqrt{m} \left( ^{(2)}\hat{\theta}_{\beta,\gamma} - \theta_2 \right) 
\rightarrow N( 0, J_{\beta}^{-1}(\theta_2) K_{\beta}(\theta_2) 
J_{\beta}^{-1}(\theta_2) ).
$$
Let $\frac{m}{m+n} \rightarrow \omega \in (0,1)$ as $m,n \rightarrow 
\infty $. Then we have
$$
\sqrt{\frac{mn}{m+n}} \left( ^{(1)}\hat{\theta}_{\beta,\gamma} - \theta_1 
\right) \rightarrow N( 0, \omega J_{\beta}^{-1}(\theta_1) 
K_{\beta}(\theta_1) J_{\beta}^{-1}(\theta_1) )
$$
and
$$
\sqrt{\frac{mn}{m+n}} \left( ^{(2)}\hat{\theta}_{\beta,\gamma} - 
\theta_2 \right) \rightarrow N( 0, (1-\omega) J_{\beta}^{-1}
(\theta_2) K_{\beta}(\theta_2) J_{\beta}^{-1}(\theta_2) ).
$$
Now, under $H_0 : \theta_1 = \theta_2$, we get that
$$
\sqrt{\frac{mn}{m+n}} \left( ^{(1)}\hat{\theta}_{\beta,\gamma} - 
^{(2)}\hat{\theta}_{\beta,\gamma} \right) \rightarrow N( 0,  
J_{\beta}^{-1}(\theta_1) K_{\beta}(\theta_1) J_{\beta}^{-1}(\theta_1) ).
$$
Next consider the second order Taylor series expansion of 
${\rm LSD}_{\beta,\gamma}\left( f_{\theta_1} , f_{\theta_2}\right)$ around 
$\theta_1 = \theta_2$ at $\left( ^{(1)}\hat{\theta}_{\beta,\gamma} , 
^{(2)}\hat{\theta}_{\beta,\gamma} \right)$ as follows
\begin{eqnarray}
    {\rm LSD}_{\beta,\gamma}\left( f_{^{(1)}\hat{\theta}_{\beta,\gamma}} , 
f_{^{(2)}\hat{\theta}_{\beta,\gamma}}\right) & = & \frac{1}{2} \sum_{i,j=1}^p 
\left( \frac{\partial^2  {\rm LSD}_{\beta,\gamma}(f_{\theta_1}, 
f_{\theta_2})}{\partial \theta_{1i} \partial \theta_{1j}}\right)_{\theta_1 
= \theta_2} \nonumber \\
& & \left(\hat{\theta}_{\beta,\gamma}^{1i} - \theta_{1i}\right) 
\left(\hat{\theta}_{\beta,\gamma}^{1j} - \theta_{1j}\right) \nonumber
\\
& & +  \sum_{i,j=1}^p \left( \frac{\partial^2  {\rm LSD}_{\beta,\gamma}
(f_{\theta_1}, f_{\theta_2})}{\partial \theta_{1i} \partial 
\theta_{2j}}\right)_{\theta_1 = \theta_2}\nonumber \\
& &  \left(\hat{\theta}_{\beta,\gamma}^{1i} 
- \theta_{1i}\right) \left(\hat{\theta}_{\beta,\gamma}^{2j} - \theta_{2j}\right) \nonumber
\\
& & +  \frac{1}{2} \sum_{i,j=1}^p \left( \frac{\partial^2  
{\rm LSD}_{\beta,\gamma}(f_{\theta_1}, f_{\theta_2})}{\partial 
\theta_{2i} \partial \theta_{2j}}\right)_{\theta_1 = \theta_2} \nonumber \\
& & 
\left(\hat{\theta}_{\beta,\gamma}^{2i} - \theta_{2i}\right) \left(\hat{\theta}_{\beta,\gamma}^{2j} - 
\theta_{2j}\right) \nonumber
\\
& & + o\left( ||^{(1)}\hat{\theta}_{\beta,\gamma} - 
\theta_1 ||^2 \right) + o\left( ||^{(2)}\hat{\theta}_{\beta,\gamma} - 
\theta_2 ||^2 \right). \nonumber 
\end{eqnarray}
But for $ i = 1, \ldots, p$, we have 
$$
\frac{\partial  {\rm LSD}_{\beta,\gamma}(f_{\theta_1}, 
f_{\theta_2})}{\partial \theta_{1i}} = \frac{1+\beta}{B} 
\left[\frac{\sum f_{\theta_{1i}}^{1+\beta}u_{\theta_{1i}}}
{\sum f_{\theta_{1i}}^{1+\beta}}
- \frac{\sum f_{\theta_{1i}}^A f_{\theta_2}^B u_{\theta_{1i}}}
{\sum f_{\theta_{1i}}^A f_{\theta_2}^B }\right].
$$
where $B = \beta - \gamma(1-\beta)$ and $A = 1+ 
\gamma(1-\beta)$ and hence 
\begin{eqnarray}
   \left( \frac{\partial^2  {\rm LSD}_{\beta,\gamma}(f_{\theta_1}, 
f_{\theta_2})}{\partial \theta_{1i} \partial \theta_{1j}}\right)_{\theta_1 = \theta_2} 
&=& 
a_{ij}^{\beta}(\theta_1), \nonumber \\
 \left( \frac{\partial^2  {\rm LSD}_{\beta,\gamma}
(f_{\theta_1}, f_{\theta_2})}{\partial \theta_{1i} 
\partial \theta_{2j}}\right)_{\theta_1 = \theta_2} 
&=& 
- a_{ij}^{\beta}(\theta_1), \nonumber \\
 \left( \frac{\partial^2  {\rm LSD}_{\beta,\gamma}(f_{\theta_1}, 
f_{\theta_2})}{\partial \theta_{2i} \partial \theta_{2j}}\right)_{\theta_1 = \theta_2} 
& =&  
a_{ij}^{\beta}(\theta_1). \nonumber
\end{eqnarray}
As in the one sample case, here also the second order partial
derivatives depend on $\beta$ only.
	Therefore, we get 
\begin{eqnarray}
    2~{\rm LSD}_{\beta,\gamma}\left( f_{^{(1)}\hat{\theta}_{\beta,\gamma}} , 
f_{^{(2)}\hat{\theta}_{\beta,\gamma}}\right) 
 &=& (^{(1)}\hat{\theta}_{\beta,\gamma} - 
\theta_1)^T A_{\beta}(\theta_1)(^{(1)}\hat{\theta}_{\beta,\gamma} - \theta_1) \nonumber \\
& & -  
2 (^{(1)}\hat{\theta}_{\beta,\gamma} - \theta_1)^T A_{\beta}(\theta_1)
(^{(2)}\hat{\theta}_{\beta,\gamma} - \theta_1) \nonumber \\
& & + (^{(2)}\hat{\theta}_{\beta,\gamma} - \theta_1)^T A_{\beta}(\theta_1)
(^{(2)}\hat{\theta}_{\beta,\gamma} - \theta_1) \nonumber \\ 
& & + o\left( ||^{(1)}
\hat{\theta}_{\beta,\gamma} - \theta_1 ||^2 \right) 
+ o\left( ||^{(2)}
\hat{\theta}_{\beta,\gamma} - \theta_2 ||^2 \right) \nonumber \\
&=&  (^{(1)}\hat{\theta}_{\beta,\gamma} - ^{(2)}\hat{\theta}_{\beta,\gamma})^T 
A_{\beta}(\theta_1)(^{(1)}\hat{\theta}_{\beta,\gamma} - 
^{(2)}\hat{\theta}_{\beta,\gamma}) \nonumber \\ 
& & + o\left( ||^{(1)}\hat{\theta}_{\beta,\gamma} 
- \theta_1 ||^2 \right) 
+ o\left( ||^{(2)}\hat{\theta}_{\beta,\gamma} \nonumber
- \theta_2 ||^2 \right), \nonumber 
\end{eqnarray}
with
$$
o\left( ||^{(1)}\hat{\theta}_{\beta,\gamma} - \theta_1 ||^2 \right) = 
o_p\left(\frac{1}{n}\right) ~~~~ \mbox{and} ~~~~  o\left( ||^{(2)}\hat{\theta}_{\beta,\gamma} - 
\theta_2 ||^2 \right) = o_p\left(\frac{1}{m}\right).
$$
Thus the asymptotic distribution of 
$$
S_{\beta,\gamma}\left( ^{(1)}\hat{\theta}_{\beta,\gamma} , ^{(2)}\hat{\theta}_{\beta,\gamma}\right) 
= \frac{2nm}{n+m} ~  {\rm LSD}_{\beta,\gamma}\left( f_{^{(1)}\hat{\theta}_{\beta,\gamma}} 
, f_{^{(2)}\hat{\theta}_{\beta,\gamma}}\right)
$$
coincides with the distribution of the random variable 
$\sum \limits_{i=1}^r ~  \zeta_i^{{\beta}}(\theta_1)Z_i^2$. 
Like the one sample case the asymptotic distribution depends 
on $\beta$ only. 
\end{proof}
\end{theorem}
We have noted that the asymptotic distribution of the LSD based test statistics 
under the simple null hypothesis is independent of the parameter $\gamma$.
Maji, Ghosh and Basu (2014) have also reported a similar observation for the robustness of 
the corresponding test statistics. They have shown that the first order influence function 
of the test statistics is always zero at the simple null and its second order influence function under null, 
being a quadratic form in the first order influence function of the minimum LSD estimator used,
is independent of the parameter $\gamma$. However, the numerical illustrations reported in 
their paper and in the next section of present paper, this independence is not true for samples
with moderate size. Therefore, as in the case of estimation, the robustness of the LSD based 
test of simple null hypothesis can not e indicated in terms of the influence function analysis
even if we even go up to second order. However, Maji, Ghosh and Basu (2014) showed 
that the robustness of the minimum LSD estimators can be measured quite accurately
in terms of the secord order influence function of the estimator. 
Extending the same idea in case of testing, it is a routine exercise to 
see that the third order influence function of the test statistics at the null, 
being a function of the second order influence function of the corresponding estimator,
can serve a better measure of robustness in this case.
In this article we have restricted ourselves 
to the simple null case. However the results may be extended 
to the case involving nuisance 
parameters following the same general approach. 
\section{A Two-Sample Example}
\label{SEC:testing:example}
Here we will discuss a two sample real data example which
is known to give rise to occasional spurious counts. This experiment
is available in  Woodruff et al.~(1984) and 
has been analyzed previously by 
Simpson (1989). 
This is a sex-linked recessive lethal experiment in drosophila (fruit flies)
to test chemical mutagenicity. Male flies were exposed either to 
$2000~\mu g$ butyraldehyde or to control conditions. The 
responses are the numbers of recessive lethal 
mutations observed among daughters of these flies. 
The data are given in Table \ref{tab:twoDdata}. 
We will use a Poisson model in this experiment
where the control responses 
are supposed to follow Poisson distribution with mean $\theta_0$
and and the treated responses follows a Poisson distribution with mean $\theta_1$. 
The two large counts for the treated group appears to be 
possible outliers. 
 We want to test
$H_0:\theta_0\geq\theta_1$ against $H_1:\theta_0<\theta_1.$  
The test statistic for testing this hypothesis is given by,
$$
^*S_{\beta,\gamma}\left( ^{(1)}\hat{\theta}_{\beta,\gamma}, 
^{(2)}\hat{\theta}_{\beta,\gamma}\right) = \frac{1}{\zeta(^{(0)}\hat{\theta}_{\beta,\gamma})} 
 \frac{2nm}{n+m} ~  {\rm LSD}_{\beta,\gamma}\left( p({^{(1)}\hat{\theta}_{\beta,\gamma}}) 
, p({^{(2)}\hat{\theta}_{\beta,\gamma}})\right),
$$
where
\begin{eqnarray}
    \zeta(^{(0)}\hat{\theta}_{\beta,\gamma}) 
&=& \frac{A_{\beta,\gamma}(^{(0)}\hat{\theta}_{\beta,\gamma}) 
K_{\beta,\gamma}(^{(0)}\hat{\theta}_{\beta,\gamma})}
{J_{\beta,\gamma}^2(^{(0)}\hat{\theta}_{\beta,\gamma})}. \nonumber \\
\end{eqnarray}
The asymptotic distribution of the statistic 
$^*S_{\gamma,\lambda}( ^{(1)}\hat{\theta}_\beta 
, ^{(2)}\hat{\theta}_\beta)$ is chi-square with one degree of freedom 
and the corresponding $p$-values are calculated and shown in the tables.
The statistic for $\beta = \gamma = 0$ case is not same as 
the likelihood ratio statistic but they are 
asymptotically same. Though the tests are different 
but the non-robust nature of the likelihood test can be seen 
under this set-up also.  
The results are shown in 
Tables \ref{tab:est2} and \ref{signed_test}. 
From the results it is evident that for $\beta \geq 0.6$ and irrespective of $\gamma$ 
the presence and the absence of outliers 
has little impact on $p$-values.
For $\gamma<0$, lower values of $\beta$ also give 
close $p$-values but for $\gamma>0$, lower values of $\beta$ the method does not perform well. 
It is clear that for large values of $\beta~(\mbox{say}~\geq 0.5)$, 
the full data and the outlier deleted data basically lead to  the same conclusion 
and almost identical $p$-values irrespective of the value of $\gamma$. 
The situation changes when $\beta$ is a small positive value close to $0$. 
In this case the role of $\gamma$ becomes decisive. Large positive values 
of $\gamma$ and small values of $\beta$ lead to a highly unstable results. 
The outlier deleted $p$-values and full data $p$-values are far from 
close in these cases. However the negative values of $\gamma$ lead to
stable inference even when $\beta=0$ or in its neighborhood. 
On the whole it appears that the two large counts in the treated 
group indicate a false significance for the likelihood ratio 
test and some other members of our class, but the 
more robust members clearly recognize the significance 
to be false.  
\begin{table}
	\centering
	\caption{Frequencies of the number of recessive lethal daughters for drosophila data}
		\begin{tabular}{c c c c c c c c c } \hline
		$x$      &0 & 1  &2  &3  &4  &5  &6 & 7 \\ \hline
Observed (Control) & 159 & 15 & 3 & 0 & 0 & 0 & 0 & 0 \\
Observed (Treated) & 110 & 11 & 5 & 0 & 0 & 0 & 1 & 1 \\
 \hline
		\end{tabular}
	 \label{tab:twoDdata}
\end{table} 
\begin{table}
	\centering
	\caption{Estimated Poisson parameters for the two-sample drosophila example;
the numbers within the parentheses 
show the corresponding estimates for the treated case after deleting 
the two outliers.
The parameter $\gamma$ is held at $0$  for the first part 
and the parameter $\beta$ is held at $0$  for the second part. 
}
	\scalebox{0.85} {	\begin{tabular}
{c c c c c c c c c c c } 
\hline
	$\beta$	        
&    $0.2$
&    $ 0.4 $   
&    $0.6 $  
&    $ 0.8 $   
& $1$ \\ \hline 
$^{(1)}\hat{\theta}_\beta$   & 0.1091 & 0.1027 & 0.099 & 0.0969 & 0.0957 \\ 
$^{(2)}\hat{\theta}_\beta$   & 0.153 & 0.1266 & 0.1143 & 0.1077 & 0.1042 \\ 
 & (0.1432) & (0.1255) & (0.1141) & (0.1105) & (0.1057) \\ 
$^{(0)}\hat{\theta}_\beta$   & 0.1264 & 0.1122 & 0.1051 & 0.1029 & 0.1 \\ 
 & (0.1229) & (0.1118) & (0.105) & (0.1028) & (0.1) \\

 \hline
	$\gamma$	        
&    $0.2$
&    $ 0.4 $   
&    $0.6 $  
&    $ 0.8 $   
& $1$ \\ \hline
$^{(1)}\hat{\theta}_\gamma$    & 0.1216 & 0.1245 & 0.1273 & 0.1303 & 0.133 \\ 
$^{(2)}\hat{\theta}_\gamma$    & 0.3227 & 0.4139 & 0.4916 & 0.5547 & 0.6059 \\ 
 & (0.1763) & (0.1854) & (0.1938) & (0.2015) & (0.2084) \\ 
$^{(0)}\hat{\theta}_\gamma$    & 0.2182 & 0.2851 & 0.3497 & 0.4055 & 0.4521 \\ 
  & (0.1444) & (0.1501) & (0.1555) & (0.1607) & (0.1655) \\ 
\hline
		\end{tabular}}
	 \label{tab:est2}
\end{table}
 \begin{table}

	\centering 
	\caption{The 
$p$-values for the two-sample drosophila data.
The outlier deleted 
$p$-values are given in the second line of 
each block
}
\scalebox{0.8}		{\begin{tabular}{c c c c c c c c c c c c} \hline
$\gamma$ \textdownarrow $\beta$ $\rightarrow$ 
& 0 & 0.1 & 0.2 & 0.3 & 0.4 & 0.5 & 0.6 & 0.7 & 0.8 & 0.9 & 1 \\ 

\hline
$-0.8$  & 0.262 & 0.341 & 0.426 & 0.509 & 0.582 & 0.642 & 0.695 & 0.738 & 0.775 & 0.806 & 0.831 \\ 
  & 0.264 & 0.343 & 0.428 & 0.511 & 0.583 & 0.643 & 0.696 & 0.739 & 0.776 & 0.807 & 0.832 \\ 
$-0.7$  & 0.235 & 0.325 & 0.417 & 0.506 & 0.582 & 0.644 & 0.698 & 0.74 & 0.776 & 0.807 & 0.831 \\ 
  & 0.239 & 0.327 & 0.42 & 0.508 & 0.584 & 0.646 & 0.7 & 0.741 & 0.777 & 0.807 & 0.832 \\ 
$-0.6$  & 0.218 & 0.314 & 0.412 & 0.501 & 0.579 & 0.645 & 0.699 & 0.742 & 0.779 & 0.807 & 0.831 \\ 
  & 0.221 & 0.317 & 0.414 & 0.503 & 0.581 & 0.647 & 0.7 & 0.744 & 0.78 & 0.808 & 0.832 \\ 
$-0.5$  & 0.207 & 0.306 & 0.406 & 0.499 & 0.578 & 0.646 & 0.7 & 0.745 & 0.779 & 0.807 & 0.831 \\ 
  & 0.211 & 0.309 & 0.408 & 0.5 & 0.58 & 0.647 & 0.701 & 0.746 & 0.78 & 0.808 & 0.832 \\ 
$-0.4$  & 0.199 & 0.301 & 0.402 & 0.496 & 0.579 & 0.647 & 0.704 & 0.746 & 0.782 & 0.809 & 0.831 \\ 
  & 0.208 & 0.305 & 0.404 & 0.497 & 0.581 & 0.649 & 0.705 & 0.747 & 0.783 & 0.81 & 0.832 \\ 
$-0.3$  & 0.187 & 0.294 & 0.398 & 0.494 & 0.578 & 0.648 & 0.702 & 0.747 & 0.782 & 0.81 & 0.831 \\ 
  & 0.211 & 0.303 & 0.403 & 0.498 & 0.58 & 0.649 & 0.704 & 0.749 & 0.783 & 0.811 & 0.832 \\ 
$-0.2$  & 0.149 & 0.284 & 0.395 & 0.493 & 0.578 & 0.647 & 0.705 & 0.749 & 0.785 & 0.81 & 0.831 \\ 
  & 0.22 & 0.307 & 0.402 & 0.496 & 0.58 & 0.648 & 0.706 & 0.75 & 0.786 & 0.811 & 0.832 \\ 
$-0.1$  & 0.058 & 0.251 & 0.388 & 0.491 & 0.576 & 0.648 & 0.705 & 0.749 & 0.785 & 0.811 & 0.831 \\ 
  & 0.235 & 0.315 & 0.404 & 0.497 & 0.579 & 0.651 & 0.706 & 0.75 & 0.786 & 0.812 & 0.832 \\ 
0 & 0.003 & 0.156 & 0.367 & 0.488 & 0.575 & 0.649 & 0.707 & 0.753 & 0.786 & 0.811 & 0.831 \\ 
   & 0.141 & 0.328 & 0.411 & 0.496 & 0.578 & 0.65 & 0.708 & 0.754 & 0.787 & 0.812 & 0.832 \\ 
0.1 & 0 & 0.03 & 0.306 & 0.479 & 0.575 & 0.648 & 0.707 & 0.753 & 0.788 & 0.814 & 0.831 \\ 
  & 0.293 & 0.345 & 0.419 & 0.5 & 0.579 & 0.65 & 0.708 & 0.754 & 0.789 & 0.814 & 0.832 \\ 
0.2 & 0 & 0.001 & 0.172 & 0.454 & 0.572 & 0.648 & 0.707 & 0.754 & 0.788 & 0.814 & 0.831 \\ 
  & 0.338 & 0.369 & 0.429 & 0.504 & 0.58 & 0.651 & 0.708 & 0.755 & 0.789 & 0.815 & 0.832 \\ 
0.3 & 0 & 0 & 0.037 & 0.391 & 0.566 & 0.648 & 0.708 & 0.754 & 0.79 & 0.815 & 0.831 \\ 
  & 0.392 & 0.398 & 0.444 & 0.51 & 0.583 & 0.65 & 0.709 & 0.755 & 0.791 & 0.815 & 0.832 \\ 
0.4 & 0 & 0 & 0.002 & 0.261 & 0.55 & 0.647 & 0.709 & 0.755 & 0.79 & 0.815 & 0.831 \\ 
  & 0.46 & 0.435 & 0.462 & 0.517 & 0.585 & 0.652 & 0.71 & 0.756 & 0.791 & 0.816 & 0.832 \\ 
0.5 & 0 & 0 & 0 & 0.101 & 0.513 & 0.644 & 0.709 & 0.757 & 0.792 & 0.816 & 0.831 \\ 
  & 0.54 & 0.478 & 0.483 & 0.526 & 0.588 & 0.652 & 0.71 & 0.758 & 0.793 & 0.817 & 0.832 \\ 
0.6 & 0 & 0 & 0 & 0.019 & 0.436 & 0.639 & 0.709 & 0.758 & 0.793 & 0.816 & 0.831 \\ 
  & 0.629 & 0.528 & 0.508 & 0.539 & 0.592 & 0.653 & 0.711 & 0.759 & 0.794 & 0.817 & 0.832 \\ 
0.7 & 0 & 0 & 0 & 0.002 & 0.3 & 0.629 & 0.709 & 0.757 & 0.793 & 0.817 & 0.831 \\ 
  & 0.727 & 0.584 & 0.538 & 0.552 & 0.597 & 0.655 & 0.712 & 0.758 & 0.794 & 0.818 & 0.832 \\ 
0.8 & 0.003 & 0 & 0 & 0 & 0.147 & 0.605 & 0.709 & 0.758 & 0.794 & 0.817 & 0.831 \\ 
  & 0.827 & 0.646 & 0.571 & 0.568 & 0.603 & 0.657 & 0.712 & 0.759 & 0.795 & 0.818 & 0.832 \\ 

\hline
		\end{tabular}}
\label{signed_test}	 
\end{table}
\section{Conclusion}
\label{Sec_Conclusion}
%
Logarithmic super divergence family acts as a super family of both LPD and LDPD family. 
The theoretical properties of this new family of divergences have been established
for discrete models and similar results under continuous set-up can be done 
in subsequent works.

\section*{Appendix}

{\bf Consistency Part:} Consider the behavior 
of ${\rm LSD}(r_n,f_{\theta})$ on a sphere $Q_a$ which has  radius 
$a$ and center at $\theta^g$. We will show, for sufficiently 
small $a$, the probability tends to one that 
$$
{\rm LSD}(r_n,f_\theta) > {\rm LSD}(r_n,f_{\theta^g}) ~~ \mbox{ for all } 
\theta ~~ \mbox{ on the surface of } Q_a
$$
so that the LSD has a local minimum with respect to $\theta$ 
in the interior of $Q_a$. At a local minimum, the estimating 
equations must be satisfied. Therefore, for any $a>0$ sufficiently 
small, the minimum LSD estimating equation have a solution $\theta_n$ 
within $Q_a$ with probability tending to one as $n \rightarrow \infty$.
	Now taking Taylor series expansion of ${\rm LSD}(r_n,f_{\theta})$ about $\theta = \theta^g$, we get 
\begin{eqnarray}
    {\rm LSD}(r_n,f_{\theta^g}) - {\rm LSD}(r_n,f_\theta)
&=& - \sum_j (\theta_j - \theta_j^g)\nabla_j {\rm LSD}(r_n,f_\theta)|_{\theta = \theta^g} \nonumber \\
& & - 
\frac{1}{2} \sum_{j,k}  (\theta_j - \theta_j^g)(\theta_k - 
\theta_k^g)\nabla_{jk} {\rm LSD}(r_n,f_\theta)|_{\theta = \theta^g} \nonumber \\  
   			& & -      \frac{1}{6} \sum_{j,k,l} (\theta_j - 
\theta_j^g)(\theta_k - \theta_k^g)(\theta_l - \theta_l^g)\nabla_{jkl} 
{\rm LSD}(r_n,f_\theta)|_{\theta = \theta^*} \nonumber \\
& & = S_1 + S_2 + S_3, ~~~~~~~~ (\mbox{say}) \nonumber
\end{eqnarray}

where $\theta^*$ lies between $\theta^g$ and $\theta$. We will now consider each terms one-by-one.

	For the Linear term $S_1$, we consider
\begin{eqnarray}
    \nabla_j {\rm LSD}(r_n,f_\theta)
= 
\frac{(1+\beta)}{A} \left[\left(\sum_x f_\theta^{1+\beta}u_{j\theta}\right)
\left(\sum_x f_\theta^{1+\beta}\right)^{-1}
-\left(\sum_x r_n^A f_\theta^{B}u_{j\theta}\right)
\left(\sum_x r_n^A f_\theta^{B}\right)^{-1}\right]. \nonumber 
\end{eqnarray}
To show  $$\nabla_j {\rm LSD}(r_n,f_\theta)
|_{\theta = \theta^g} \xrightarrow{P} 0,$$ we need to show 
  $$\sum_x K(\delta_n^g(x)) f_{\theta^g}^{1+\beta}(x) u_{j\theta^g}(x)
\xrightarrow{P} \sum_x K(\delta_g^g(x)) f_{\theta^g}^{1+\beta}(x) u_{j\theta^g}(x)$$
and 
  $$\sum_x K(\delta_n^g(x)) f_{\theta^g}^{1+\beta}(x) 
\xrightarrow{P} \sum_x K(\delta_g^g(x)) f_{\theta^g}^{1+\beta}(x) $$
as $n\rightarrow\infty$ and 
where $\delta_n^g(x)$ is the $\delta_n(x)$ evaluated at $\theta = \theta^g$. We will now 
show that 
\begin{eqnarray}
\label{EQ:k_delta_converge}
    \sum_x K(\delta_n^g(x))f_{\theta^g}^{1+\beta}(x)
u_{j\theta^g}(x) \xrightarrow{P} \sum_x K(\delta_g^g(x)) 
f_{\theta^g}^{1+\beta}(x) u_{j\theta^g}(x)
\end{eqnarray}
as $n \rightarrow \infty$. 
Note that by assumption $7$ and the fact that $r_n(x) \rightarrow g(x) $ 
almost surely (a.s.) by Strong law of large numbers (SLLN), it follows that 
\begin{eqnarray}
\label{EQ:k_prime_bound}
	|K'(\delta)| = |A||\delta|^{A-1} < 2|A|C = C_1,  ~~~~~~ \mbox{(say)}
\end{eqnarray}
for any $\delta$ in between $\delta_n^g(x)$ and $\delta_g^g(x)$ 
(uniformly in $x$). So,  by using the one-term Taylor series expansion,
\begin{eqnarray}
 & &   |\sum_x K(\delta_n^g(x))f_{\theta^g}^{1+\beta}(x)u_{j\theta^g}(x) - 
\sum_x K(\delta_g^g(x)) f_{\theta^g}^{1+\beta}(x) u_{j\theta^g}(x) |  \nonumber \\
&\le & C_1 \sum_x |\delta_n^g(x) - \delta_g^g(x)| 
f_{\theta^g}^{1+\beta}(x) |u_{j\theta^g}(x) |. \nonumber \\ \nonumber 
\end{eqnarray}

	However, by Lemma \ref{lem:lamma1}(1), 
\begin{eqnarray}
     E\left[|\delta_n^g(x) - \delta_g^g(x)|\right] \le 
\frac{(g(x)(1-g(x))^{1/2}}{f_{\theta^g}(x) 
\sqrt n} \rightarrow 0  ~~~ as ~ ~ n \rightarrow \infty.
\end{eqnarray} 
	and, by Lemma \ref{lem:lamma1}(2), we have 
\begin{eqnarray}
 & &  E\left[ C_1 \sum_x |\delta_n^g(x) - \delta_g^g(x)| 
f_{\theta^g}^{1+\beta}(x) |u_{j\theta^g}(x) |\right] \nonumber \\
&\le& 2C_1 \sum_x g^{1/2}(x) f_{\theta^g}^{\beta}(x) |u_{j\theta^g}(x) |  ~ < ~ \infty     
 ~~~~~~~~~\mbox{(by Assumption $5$)}.
\end{eqnarray}
	Hence, by DCT, we get, 
\begin{eqnarray}
     E\left[|\sum_x K(\delta_n^g(x))f_{\theta^g}^{1+\beta}(x)
u_{j\theta^g}(x) - \sum_x K(\delta_g^g(x)) f_{\theta^g}^{1+\beta}(x) 
u_{j\theta^g}(x) |\right] \rightarrow 0 ~~~~ as ~ n \rightarrow \infty \nonumber
\end{eqnarray}
so that by Markov inequality we have the desired claim. 

By similar argument we can show 
\begin{eqnarray}
 \sum_x K(\delta_n^g(x)) f_{\theta^g}^{1+\beta}(x) \nonumber \\
\xrightarrow{P} \sum_x K(\delta_g^g(x)) f_{\theta^g}^{1+\beta}(x). 
\end{eqnarray}
Therefore, we have
\begin{eqnarray}
    \nabla_j {\rm LSD}(r_n,f_\theta)|_{\theta = \theta^g} \xrightarrow{P} 0.
\end{eqnarray}
	Thus, with probability tending to one, $|S_1| < p a^3$, where 
$p$ is the dimension of $\theta$ and $a$ is the radius of $Q_a$.

	Next we consider the quadratic term $S_2$. We have,
\begin{eqnarray}
\label{EQ:S2_form}
\frac{A}{1+\beta}\nabla_{jk} {\rm LSD}(r_n, f_\theta) &=& 
(1+\beta)\left(\sum_x f_\theta^{1+\beta} u_{j\theta}u_{k\theta}\right) 
\left(\sum_x f_\theta^{1+\beta} \right)^{-1} \nonumber \\ 
& & - (1+\beta)\left(\sum_x f_\theta^{1+\beta} u_{j\theta}\right)
\left(\sum_x f_\theta^{1+\beta} u_{k\theta}\right) 
\left(\sum_x f_\theta^{1+\beta} \right)^{-2}  \nonumber \\ 
& & -B \left(\sum_x r_n^A f_\theta^B u_{j\theta} u_{k\theta}\right)
\left(\sum_x r_n^A f_\theta^B\right)^{-1} \nonumber \\ 
& & +B \left( \sum_x r_n^A f_\theta^B u_{j\theta}\right)\left(\sum_x r_n^A f_\theta^B u_{k\theta}\right)
\left(\sum_x r_n^A f_\theta^B \right)^{-2} 
\end{eqnarray}
We will now show that 
\begin{eqnarray}
\nabla_{jk} {\rm LSD}(r_n, f_\theta) |_{\theta=\theta^g} \xrightarrow{P} - J_g^{j,k}. 
\end{eqnarray}
To show this we will show that 
\begin{eqnarray}
\label{EQ:k_ujk}
  \sum_x K(\delta_n^g(x)) f_{\theta^g}^{1+\beta}(x) u_{j\theta^g}(x) u_{k\theta^g}(x) \nonumber \\
\rightarrow \sum_x K(\delta_g^g(x)) f_{\theta^g}^{1+\beta}(x) u_{j\theta^g}(x) u_{k\theta^g}(x).
\end{eqnarray}
Note that by assumption $7$ and using the fact that $r_n(x) \rightarrow g(x) ~~ a.s. $ by SLLN
and by (\ref{EQ:k_prime_bound}) we get the following expansion,
\begin{eqnarray}
 & &   |\sum_x K(\delta_n^g(x))f_{\theta^g}^{1+\beta}(x)u_{j\theta^g}(x)u_{k\theta^g}(x)  - 
\sum_x K(\delta_g^g(x)) f_{\theta^g}^{1+\beta}(x) u_{j\theta^g}(x)u_{k\theta^g}(x)|  \nonumber \\
&\le & C_1 \sum_x |\delta_n^g(x) - \delta_g^g(x)| 
f_{\theta^g}^{1+\beta}(x) |u_{j\theta^g}(x)u_{k\theta^g}(x)|. \nonumber \\ \nonumber 
\end{eqnarray}

	However, by Lemma \ref{lem:lamma1}(1), 
\begin{eqnarray}
     E\left[|\delta_n^g(x) - \delta_g^g(x)|\right] \le 
\frac{(g(x)(1-g(x))^{1/2}}{f_{\theta^g}(x) 
\sqrt n} \rightarrow 0  ~~~ as ~ ~ n \rightarrow \infty.
\end{eqnarray} 
	and, by Lemma \ref{lem:lamma1}(2), we have 
\begin{eqnarray}
 & &  E\left[ C_1 \sum_x |\delta_n^g(x) - \delta_g^g(x)| 
f_{\theta^g}^{1+\beta}(x) |u_{j\theta^g}(x) u_{k\theta^g}(x)|\right] \nonumber \\
&\le& 2C_1 \sum_x g^{1/2}(x) f_{\theta^g}^{\beta}(x) |u_{j\theta^g}(x) u_{k\theta^g}(x)|  ~ < ~ \infty     
 ~~~~~~~~~\mbox{(by Assumption $5$)}. \nonumber
\end{eqnarray}
	Hence, by DCT, we get, 
\begin{eqnarray}
\label{EQ:s2_first}
|\sum_x K(\delta_n^g(x))f_{\theta^g}^{1+\beta}(x)
u_{j\theta^g}(x) u_{k\theta^g}(x) \nonumber  \\  - \sum_x K(\delta_g^g(x)) f_{\theta^g}^{1+\beta}(x) 
u_{j\theta^g}(x) u_{k\theta^g}(x)|
& & \rightarrow 0 
\end{eqnarray}
 as $n \rightarrow \infty$.
so that by Markov inequality we have the desired claim. 
and along with we will use 
\begin{eqnarray}
\label{EQ:s2_second}
 \sum_x K(\delta_n^g(x)) f_{\theta^g}^{1+\beta}(x) u_{j\theta^g}(x)
\xrightarrow{P} \sum_x K(\delta_g^g(x)) f_{\theta^g}^{1+\beta}(x) u_{j\theta^g}(x)
\end{eqnarray}
and 
\begin{eqnarray}
\label{EQ:s2_third}
\sum_x K(\delta_n^g(x)) f_{\theta^g}^{1+\beta}(x) 
\xrightarrow{P} \sum_x K(\delta_g^g(x)) f_{\theta^g}^{1+\beta}(x) 
\end{eqnarray}
as $n\rightarrow\infty$ from the previous part of the proof. 
	Thus, combining 
(\ref{EQ:s2_first}), (\ref{EQ:s2_second}) and (\ref{EQ:s2_third}), 
we get that,
\begin{eqnarray}
\nabla_{jk} {\rm LSD}(r_n, f_\theta) |_{\theta=\theta^g} \xrightarrow{P} - J_g^{j,k}.
\end{eqnarray}
 Therefore,
\begin{eqnarray}
\label{EQ:2s2_form}
    2S_2 &=& \frac{1+\beta}{A} \sum_{j,k} 
\left\{ \nabla_{jk} {\rm LSD}(r_n, f_\theta) |_{\theta=\theta^g} 
 - 
( - J_g^{j,k})\right\} (\theta_j - \theta_j^g)(\theta_k - \theta_k^g) \nonumber \\
 & &  ~~~~~~~ + \sum_{j,k} \left\{ - \left(\frac{(1+\beta)}{A} 
J_g^{j,k} \right) (\theta_j - \theta_j^g)(\theta_k - \theta_k^g) \right\}. \\ \nonumber
\end{eqnarray}
	Now the absolute value of the first term in above (\ref{EQ:2s2_form}) is 
$ < p^2 a^3$ with probability tending to one. And, the second term in 
(\ref{EQ:2s2_form}) is a negative definite quadratic form in the variables 
$(\theta_j - \theta_j^g)$. Letting $\gamma_1$ be the largest 
eigenvalue of $\frac{(1+\beta)}{A} J_g$, the quadratic form is 
$ < \gamma_1 a^2$. Combining the two terms, we see that there 
exists $c > 0$ and $a_0 > 0$ such that for $a < a_0$, we have  
$ S_2 < -c a^2$ with probability tending to one.
	Finally, considering the cubic term $S_3$, we have
\begin{align*}
 \frac{A}{1+\beta} \nabla_{jkl} {\rm LSD}(r_n, f_\theta)
= (1+\beta)^2 \left(\sum_x f_\theta^{1+\beta} u_{j\theta}u_{k\theta}u_{l\theta}\right)
\left(\sum_x f_\theta^{1+\beta}\right)^{-1} \nonumber \\
+(1+\beta)\left(\sum_x f_\theta^{1+\beta}
u_{jl\theta}u_{k\theta}\right) 
\left(\sum_x f_\theta^{1+\beta}\right)^{-1}
+(1+\beta)\left(\sum_x f_\theta^{1+\beta}
u_{j\theta}u_{kl\theta}\right) 
\left(\sum_x f_\theta^{1+\beta}\right)^{-1} \nonumber \\
-(1+\beta)^2 \left(\sum_x f_\theta^{1+\beta}u_{j\theta}u_{k\theta}\right)
\left(\sum_x f_\theta^{1+\beta}u_{l\theta}\right)\left(\sum_x f_\theta^{1+\beta}\right)^{-2}
+(1+\beta)\left(\sum_x f_\theta^{1+\beta}u_{l\theta}u_{jk\theta}\right)
\left(\sum_x f_\theta^{1+\beta}\right)^{-1}\nonumber \\
+\left(\sum_x f_\theta^{1+\beta}u_{jkl\theta}\right)\left(\sum_x f_\theta^{1+\beta}\right)^{-1}
-(1+\beta)\left(\sum_x f_\theta^{1+\beta}u_{jk\theta}\right)
\left(\sum_x f_\theta^{1+\beta}u_{l\theta}\right)\left(\sum_x f_\theta^{1+\beta}\right)^{-2}\nonumber \\
-(1+\beta)^2\left(\sum_x f_\theta^{1+\beta}u_{j\theta}u_{l\theta}\right)
\left(\sum_x f_\theta^{1+\beta}u_{k\theta}\right)\left(\sum_x f_\theta^{1+\beta}\right)^{-2}\nonumber \\
-(1+\beta)\left(\sum_x f_\theta^{1+\beta}u_{jl\theta}\right)\left(\sum_x f_\theta^{1+\beta}u_{k\theta}\right)
\left(\sum_x f_\theta^{1+\beta}\right)^{-2}\nonumber \\
-(1+\beta)^2\left(\sum_x f_\theta^{1+\beta}u_{j\theta}\right)\left(\sum_x f_\theta^{1+\beta}u_{k\theta}u_{l\theta}\right)
\left(\sum_x f_\theta^{1+\beta}\right)^{-2}\nonumber \\
-(1+\beta)\left(\sum_x f_\theta^{1+\beta}u_{j\theta}\right)\left(\sum_x f_\theta^{1+\beta}u_{kl\theta}\right)
\left(\sum_x f_\theta^{1+\beta}\right)^{-2}\nonumber \\
+2(1+\beta)^2\left(\sum_x f_\theta^{1+\beta}u_{j\theta}\right)\left(\sum_x f_\theta^{1+\beta}u_{k\theta}\right)
\left(\sum_x f_\theta^{1+\beta}u_{l\theta}\right)\left(\sum_x f_\theta^{1+\beta}\right)^{-3}\nonumber \\
-B^2\left(\sum_x r_n^A f_\theta^B u_{j\theta}u_{k\theta}u_{l\theta}\right)
\left(\sum_x r_n^A f_\theta^B\right)^{-1}
-B\left(\sum_x r_n^A f_\theta^B u_{jl\theta}u_{k\theta}\right)\left(\sum_x r_n^A f_\theta^B\right)^{-1}\nonumber \\
-B\left(\sum_x r_n^A f_\theta^B u_{j\theta}u_{kl\theta}\right)\left(\sum_x r_n^A f_\theta^B\right)^{-1}
+B^2\left(\sum_x r_n^Af_\theta^B u_{j\theta}u_{k\theta}\right)\left(\sum_x r_n^Af_\theta^B u_{l\theta}\right)
\left(\sum_x r_n^A f_\theta^B\right)^{-2}\nonumber \\
-B\left(\sum_x r_n^Af_\theta^B u_{jk\theta}u_{l\theta}\right)\left(\sum_x r_n^Af_\theta^B u_{jkl\theta}\right)
\left(\sum_x r_n^A f_\theta^B\right)^{-2}\nonumber \\
+B\left(\sum_x r_n^Af_\theta^B u_{jk\theta}\right)\left(\sum_x r_n^Af_\theta^B u_{l\theta}\right)
\left(\sum_x r_n^A f_\theta^B\right)^{-2}\nonumber \\
+B^2\left(\sum_x r_n^Af_\theta^B u_{j\theta}u_{l\theta}\right)\left(\sum_x r_n^Af_\theta^B u_{k\theta}\right)
\left(\sum_x r_n^A f_\theta^B\right)^{-2}\nonumber \\
+B\left(\sum_x r_n^Af_\theta^B u_{jl\theta}\right)\left(\sum_x r_n^Af_\theta^B u_{k\theta}\right)
\left(\sum_x r_n^A f_\theta^B\right)^{-2}\nonumber \\
+B^2\left(\sum_x r_n^Af_\theta^B u_{k\theta}u_{l\theta}\right)\left(\sum_x r_n^Af_\theta^B u_{j\theta}\right)
\left(\sum_x r_n^A f_\theta^B\right)^{-2}\nonumber \\\displaybreak[3]\nonumber \\
+B\left(\sum_x r_n^Af_\theta^B u_{kl\theta}\right)\left(\sum_x r_n^Af_\theta^B u_{j\theta}\right)
\left(\sum_x r_n^A f_\theta^B\right)^{-2}\nonumber \\
-2B^2\left(\sum_x r_n^Af_\theta^B u_{j\theta}\right)\left(\sum_x r_n^Af_\theta^B u_{k\theta}\right)
\left(\sum_x r_n^Af_\theta^B u_{l\theta}\right)
\left(\sum_x r_n^A f_\theta^B\right)^{-3}.\nonumber \\
\end{align*}
Using the assumptions and $r_n(x)\rightarrow g(x)$ we can show the 
cubic term $S_3$ is also bounded. 
Hence, we have $|S_3|< b a^3$ on the sphere $Q_a$ with probability tending to one. 
	Combining the three inequality we get that 
$$
max(S_1 + S_2 + S_3 ) < - c a^2 + (b + p) a^3   
~~\left[~ < 0 ~~\mbox{  for } ~~ a< \frac{c}{b+p} \right].
$$
Thus, for any sufficiently small $a$, there 
exists a sequence of roots $\theta_n = \theta_n(a)$ to 
the minimum LSD estimating equation such that 
$P(||\theta_n - \theta^g||_2 < a)$ converges to one, 
where $||.||_2$ denotes the $L_2-$norm.
	It remains to show that we can determine such a 
sequence independent of $a$. For let $\theta_n^*$ be the 
root which is closest to $\theta^g$. This exists because 
the limit of a sequence of roots is again a root by the 
continuity of the LSD. Hence proved the consistency part.
%
%

\begin{thebibliography}{99}
\bibitem{P00}
Pearson, K. (1900). On the criterion that a given system of deviations from the probable in
the case of a correlated system of variables is such that it can be reasonably supposed to
have arisen from random sampling. Philosophical Magazine, {\bf 50}, 157--175.
\end{thebibliography}

\end{document}